\documentclass[12pt,draft]{amsart}

\makeatletter

\theoremstyle{plain}    
\newtheorem{thm}{Theorem}[section]
\numberwithin{equation}{section} 
\numberwithin{figure}{section} 
\theoremstyle{plain}    
\newtheorem*{thm*}{Theorem} 
\theoremstyle{plain}    
\theoremstyle{plain}    
\newtheorem{lem}[thm]{Lemma} 
\theoremstyle{plain}    
\theoremstyle{definition}
\newtheorem{defn}[thm]{Definition}
\theoremstyle{remark}
\newtheorem{rem}[thm]{Remark}
\theoremstyle{remark}

\theoremstyle{remark}    

\theoremstyle{remark}    

\theoremstyle{definition}  

\theoremstyle{remark}
  \newtheorem*{acknowledgement*}{Acknowledgement} 

\theoremstyle{plain}    

\theoremstyle{plain}    
\theoremstyle{plain}    
\theoremstyle{plain}    
\theoremstyle{definition}

\theoremstyle{remark}

\theoremstyle{remark}    

\theoremstyle{remark}    

\theoremstyle{plain}    

\usepackage{amssymb} \usepackage{amscd}

\makeatother

\begin{document}

\title[AF embeddings]{AF embeddings and the numerical
computation of spectra in irrational rotation algebras}

\author{Nathanial P. Brown}

\address{Department of Mathematics, Penn State University, State
College, PA 16802}

\email{nbrown@math.psu.edu}

\thanks{Partially supported by DMS-0244807.}


\begin{abstract} 
The spectral analysis of discretized one-dimensional Schr\"{o}dinger
operators is a very difficult problem which has been studied
by numerous mathematicians.  A natural problem at the interface of
numerical analysis and operator theory is that of finding finite
dimensional matrices whose eigenvalues approximate the spectrum of an
infinite dimensional operator.  In this note we observe that the
seminal work of Pimsner-Voiculescu on AF embeddings of irrational
rotation algebras provides a nice answer to the finite dimensional
spectral approximation problem for a broad class of operators
including the quasiperiodic case of the Schr\"{o}dinger operators
mentioned above.  Indeed, the theory of continued fractions not only
provides good matrix models for spectral computations (i.e.\ the
Pimsner-Voiculescu construction) but also yields {\em sharp} rates of
convergence for spectral approximations of operators in irrational
rotation algebras.
\end{abstract}

\maketitle

\section{Introduction}

In this paper we address the problem of finding numerical
approximations to the spectrum of a bounded linear operator on a
complex, separable Hilbert space.  This is a very difficult problem in
general and there are various ways to attack it as well as various
quantities that one may wish to use to approximate the spectrum (e.g.\
pseudospectra or numerical ranges).  A very natural approach is to
start with a given operator $T$ acting on a Hilbert space $H$ and
``compress'' it to a finite dimensional subspace of $H$.  This
compression is a finite dimensional matrix whose eigenvalues can
hopefully be numerically computed and one further hopes that these
eigenvalues will somehow approximate the spectrum of the infinite
dimensional operator $T$.  This approach has been studied by numerous
authors (see, for example, \cite{HRS}, \cite{bottcher} and the
references therein) and a number of interesting results have been
obtained for large classes of operators. 

Though the present work grew naturally out of the author's own study
of the method described above (cf.\ \cite[Section
6]{brown:finitesection}) it turns out that an absurdly abstract
approach to numerical approximations of spectra actually proves useful
in at least one important case.  The main result of this paper simply
observes how some very specialized and technical work in C$^*$-algebra
theory can be used to get excellent approximations to spectra of
some important operators.  In particular, this
applies to discretized
one-dimensional Schr\"{o}dinger operators with quasiperiodic
potential. See \cite{boca:book} for a nice treatment of the spectral
theory of Almost Mathieu operators (up to 2001) and
\cite{bourgain-jitomirskaya} (and its references) for a recent
contribution to the more general case of quasiperiodic Schr\"{o}dinger
operators.

When presented abstractly the basic idea of this note becomes quite
simple, but will require a bit of C$^*$-algebra theory to explain.
The first basic fact we need is that all finite dimensional
C$^*$-algebras have a special form -- they are just finite direct sums
of finite dimensional matrix algebras.  We will also need the 
definition of an AF (``approximately finite dimensional'') algebra; A
C$^*$-algebra $C$ is AF if there exist finite dimensional
C$^*$-subalgebras $C_1 \subset C_2 \subset \cdots \subset C$ whose
union is dense in $C$.  The observation below is the trivial, but key,
idea of this paper. 

{\noindent\bf Observation}: If $T \in B(H)$ is an element in an AF
C$^*$-subalgebra of $B(H)$ then there exist finite dimensional
matrices whose spectral quantities (e.g.\ pseudospectra or numerical
ranges) approximate those of $T$.

Indeed, if $T \in C \subset B(H)$ with $C$ an AF algebra then we can
find operators $T_n$ such that $\|T - T_n \| \to 0$ (hence spectral
quantities of $T_n$ are close to those of $T$) and each $T_n$ is
contained in a finite dimensional C$^*$-algebra (hence can be
identified with a finite dimensional matrix).  There is one subtlety
that must be mentioned.  Namely, the actual spectrum of $T_n$ need not
be close to the actual spectrum of $T$ (unless $T$ happens to be
normal) but singular values, numerical ranges and pseudospectra will
be close.  Moreover, the actual spectrum of $T$ is equal to the
intersection of all its pseudospectra and hence we still get
reasonable approximations to the spectrum of $T$ by looking at
pseudospectra of the $T_n$'s.  In the case both $T$ and
$T_n$ are normal (e.g.\ self-adjoint Schr\"{o}dinger operators) then
the actual spectra are close in a very strong sense; the Hausdorff
distance between their spectra is bounded by $\| T - T_n \|$.

AF-embeddability (i.e.\ studying which operators belong to AF
algebras) has been considered by numerous authors but most of this
work is too general to be of much help in finding {\em explicit}
matrices which approximate given operators.  For example, it is an
easy consequence of the spectral theorem that every normal operator is
contained in an AF-subalgebra of $B(H)$ but this gives no clue how to
actually find the right matrix approximations. However, returning to
the seminal paper on AF-embeddability (cf.\ \cite{PV}) one finds that
not only are explicit matrix models provided (for operators in
irrational rotation algebras) but computable and {\em sharp} rates of
convergence can be proved.

In section 2 we will review the necessary aspects of the theory of
continued fractions.  We also state a few simple lemmas that will be
needed later. 

In section 3 we review the construction of Pimsner-Voiculescu which
provides explicit matrix models for AF-embeddings of the irrational
rotation algebras.  We will not reproduce the proofs of any estimates
as these can be found either in the original paper \cite{PV} or in
\cite[Chapter VI]{davidson}.  

In section 4 we state and prove the main result for spectral
approximations of operators in irrational rotation algebras.  The
proof really amounts to working out a few estimates as the hard part
was accomplished more than 20 years ago in \cite{PV}. 

Finally, in section 5 we show how results of Haagerup-R{\o}rdam can be
used to give general ``one-sided'' rates of convergence.  Though the
results of this section are not as good as previous sections they have
the advantage of always being numerically implementable.

\section{Preliminaries}

In this section we will review some basic facts about continued
fractions and point out a few technical (but very elementary)
inequalities that we will need later.  Though there are a number of
excellent books on continued fractions, a classic reference is
\cite{hardy-wright}.  The facts we will need, however, can be found in
any book on the subject.  We also state the definition of
pseudospectrum at the very end of this section.

Given an irrational number $\theta \in (0,1)$ there exist unique
positive integers $a_1, a_2, \ldots$ such that $$\theta = \lim_{n\to
\infty} [a_1, \ldots, a_n] = \lim_{n\to \infty} \frac{1}{a_1 +
\frac{1}{a_2 + \frac{1}{\ddots + \frac{1}{a_n}}}}.$$ Note that if
$\theta$ is given then one can compute the integers $a_n$ via the
following recursive procedure: Define $\theta_1 = \frac{1}{\theta}$,
$a_1 = \lfloor \theta_1 \rfloor$, $\theta_{n+1} = \frac{1}{\theta_n -
a_n}$, $a_{n+1} = \lfloor \theta_{n+1} \rfloor$, where $\lfloor x
\rfloor$ is the integer part of a real number $x$.  As is customary,
we use $p_n$, $q_n$ to denote the numerator and denominator,
respectively, of the $n^{th}$ {\em convergent}: $$[a_1, \ldots, a_n] =
\frac{1}{a_1 + \frac{1}{a_2 + \frac{1}{\ddots + \frac{1}{a_n}}}} =
\frac{p_n}{q_n}.$$ The $p_n$'s and $q_n$'s satisfy the recursive
relations $$p_n = a_n p_{n-1} + p_{n-2}, \ \ q_n = a_n q_{n-1} +
q_{n-2},$$ where $p_0 = 0$, $p_1 = 1$, $q_0 = 1$ and $q_1 = a_1$.  The
following remarkable fact can be found in any book on continued
fractions. 

\begin{thm} 
\label{thm:contfracestimate}
For every $n$ one has $| \theta - \frac{p_n}{q_n} | <
\frac{1}{q_n q_{n+1}} < (\frac{1}{q_n})^2$.
\end{thm}

We now record a few trivial lemmas that will be
used in later estimates.  

\begin{lem} If $F(k)$ is the $k^{th}$ Fibonacci number (where $F(0) = 1
= F(1)$, $F(2) = 2$, etc.) then for each pair of positive integers
$n$ and $k$ we have $q_{n+k} \geq F(k)q_n$.
\end{lem}

\begin{proof} Use induction and the recursion formula.
\end{proof}

\begin{lem} $\sum\limits_{k = 0}^{\infty} \frac{1}{F(k)} \leq
\frac{2\sqrt{5}}{\sqrt{5} - 1}$.
\end{lem}

\begin{proof} Induction shows that $F(k) \geq \big(\frac{1 +
\sqrt{5}}{2})^{k-1}$ and hence the result follows by taking
reciprocals and applying the formula for a geometric series.
\end{proof}

\begin{lem} 
\label{thm:lemma1}
$\sum\limits_{k = 0}^{\infty} \frac{1}{q_{n+k}} \leq
\frac{1}{q_n} \frac{2\sqrt{5}}{\sqrt{5} - 1}$.
\end{lem}

\begin{proof} Apply the two previous lemmas.
\end{proof}



We close this section with the definition of pseudospectra.

\begin{defn}
\label{thm:pseudo}
For $\varepsilon > 0$, the $\varepsilon$-{\em pseudospectrum} of $T
\in B(H)$ is $$\sigma^{(\varepsilon)}(T) = \{ \lambda \in {\mathbb C}:
\| (\lambda - T)^{-1} \| \geq \frac{1}{\varepsilon} \},$$ where $\|
X^{-1}\| = \infty$ if $X \in B(H)$ is a non-invertible operator. Note
that the usual spectrum $\sigma(T)$ is contained in
$\sigma^{(\varepsilon)}(T)$ for every $\varepsilon > 0$ and, moreover,
$$\sigma(T) = \bigcap_{\varepsilon > 0} \sigma^{(\varepsilon)}(T).$$
See {\bf http://web.comlab.ox.ac.uk/projects/pseudospectra} for more,
including software for computing pseudospectra.
\end{defn}


\section{The Pimsner-Voiculescu Construction}

In this section we state the technical aspects of the
Pimsner-Voiculescu construction that we will need.  We see no reason to
reproduce proofs as all the details can be found in \cite[Sections
VI.4 and VI.5]{davidson}.

Throughout this note, $M_k ({\mathbb C})$ will denote the $k\times k$
complex matrices.  Let $0 < p < q$ be integers, $\omega = e^{2\pi i
\frac{p}{q}}$ and define matrices $u_{\frac{p}{q}}, v_{\frac{p}{q}}
\in M_{q}({\mathbb C})$ as follows.
$$u_{\frac{p}{q}} = 
\begin{pmatrix}
0 & 1 & 0 & 0 & \cdots & 0 & 0\\
0 & 0 & 1 & 0 & \cdots & 0 & 0\\
0 & 0 & 0 & 1 & \cdots & 0 & 0\\
0 & 0 & 0 & 0 & \cdots & 0 & 0\\
\vdots & \vdots & \vdots & \vdots & \ddots & \vdots & \vdots\\ 
0 & 0 & 0 & 0 & \cdots & 0 & 1\\
1 & 0 & 0 & 0 & \cdots & 0 & 0\\
\end{pmatrix}$$
$$v_{\frac{p}{q}} = 
\begin{pmatrix}
1 & 0 & 0 & 0 & \cdots & 0 & 0\\
0 & \omega & 0 & 0 & \cdots & 0 & 0\\
0 & 0 & \omega^{2} & 0 & \cdots & 0 & 0\\
0 & 0 & 0 & \omega^{3} & \cdots & 0 & 0\\
\vdots & \vdots & \vdots & \vdots & \ddots & \vdots & \vdots\\ 
0 & 0 & 0 & 0 & \cdots & \omega^{q-2} & 0\\
0 & 0 & 0 & 0 & \cdots & 0 & \omega^{q-1}\\
\end{pmatrix}$$

If we are given an irrational number $\theta \in (0,1)$ and we let
$\frac{p_n}{q_n}$ denote the $n^{th}$ convergent (as in the previous
section) then to ease notation we will let $\theta_n =
\frac{p_n}{q_n}$, $u_{\theta_n} = u_{\frac{p_n}{q_n}}$ and
$v_{\theta_n} = v_{\frac{p_n}{q_n}}$.  

\begin{thm}(Pimsner-Voiculescu) For each natural number $n$ there is a  
$*$-monomorphism $\varphi_n : M_{q_{n-1}}({\mathbb C}) \oplus
M_{q_n}({\mathbb C}) \hookrightarrow M_{q_n}({\mathbb C}) \oplus
M_{q_{n+1}}({\mathbb C})$ such that $$\| \varphi_n (u_{\theta_{n-1}}
\oplus u_{\theta_n}) - u_{\theta_n} \oplus u_{\theta_{n+1}} \| <
\frac{2\pi}{q_{n-1}}$$ and $$\| \varphi_n (v_{\theta_{n-1}} \oplus
v_{\theta_n}) - v_{\theta_n} \oplus v_{\theta_{n+1}} \| <
\frac{\pi}{q_{n-1}} + \frac{4\pi}{q_n}.$$
\end{thm}

\section{Computing Spectra in Irrational Rotation Algebras}

For an arbitrary number $\theta \in (0,1)$ there is a universal
C$^*$-algebra, denoted by $A_{\theta}$, which is generated by two
unitaries $U_{\theta}, V_{\theta}$ subject to the commutation relation
$$U_{\theta}V_{\theta} = e^{2\pi i \theta}V_{\theta}U_{\theta}.$$
Universality means that if $\tilde{U}_{\theta}, \ \tilde{V}_{\theta}$
are any other unitary operators satisfying the same relation then
there exists a surjective $*$-homomorphism $A_{\theta} \to
C^*(\tilde{U}_{\theta}, \ \tilde{V}_{\theta})$ such that $U_{\theta}
\mapsto \tilde{U}_{\theta}$ and $V_{\theta} \mapsto
\tilde{V}_{\theta}$.  One crucial fact which we will need is that if
$\theta$ is irrational then $A_{\theta}$ is a simple C$^*$-algebra
(i.e.\ has no non-trivial, closed, two-sided ideals) and hence {\em
any} $*$-homomorphism will be injective in this case (cf.\
\cite[Theorem 1.10]{boca:book}).

In the theorem below $d_{H} (\cdot, \cdot)$ denotes the Hausdorff
distance between two compact subsets of the complex plane: $$d_H
(\Sigma, \Lambda) = \max \{ \sup_{\sigma \in \Sigma} d(\sigma,
\Lambda), \sup_{\lambda \in \Lambda} d(\lambda, \Sigma) \},$$ where
$d(\sigma, \Lambda) = \inf_{\lambda \in \Lambda} |\sigma - \lambda|$. 

\begin{thm} Let $\theta \in (0,1)$ be irrational, $\theta_n =
\frac{p_n}{q_n}$, $\alpha_{\pm 1}, \beta_{\pm 1} \in {\mathbb C}$ be
four complex numbers and $M = \max\{|\alpha_{\pm 1}|, |\beta_{\pm 1}|\}$.
If we let $H_{\theta} = \alpha_1 U_{\theta} + \alpha_{-1} U_{\theta}^*
+ \beta_1 V_{\theta} + \beta_{-1} V_{\theta}^*$ and $h_{\theta_n} =
\alpha_1 u_{\theta_{n}} + \alpha_{-1} u_{\theta_{n}}^* + \beta_1
v_{\theta_{n}} + \beta_{-1} v_{\theta_n}^*$ then:
\begin{enumerate}
\item (cf.\ Definition \ref{thm:pseudo}) For each
$\varepsilon > 0$ we have
$$\bigg(\sigma^{(\varepsilon)}(h_{\theta_{n-1}}) \cup
\sigma^{(\varepsilon)}(h_{\theta_{n}})\bigg) \subset \sigma^{(\varepsilon +
\varepsilon_n)}(H_{\theta}) \subset \bigg(\sigma^{(\varepsilon +
2\varepsilon_n)}(h_{\theta_{n-1}}) \cup \sigma^{(\varepsilon +
2\varepsilon_n)}(h_{\theta_{n}})\bigg),$$ where $\varepsilon_n =
204M(\frac{1}{q_{n-1}} + \frac{1}{q_{n}})$.


\item (Normal Case) If $h_{\theta_{n-1}}, h_{\theta_{n}}$ and
$H_{\theta}$ all happen to be normal operators (e.g.\ self-adjoints)
then $$d_H (\sigma (h_{\theta_{n-1}}) \cup \sigma (h_{\theta_{n}}) ,
\sigma (H_{\theta})) \leq 204M(\frac{1}{q_{n-1}} + \frac{1}{q_{n}}).$$
\end{enumerate} 
\end{thm}

\begin{proof}  It follows from the Pimsner-Voiculescu 
construction that there is an AF algebra ${\mathfrak A}_{\theta}$ and
a sequence of $*$-monomorphisms $\psi_n : M_{q_{n-1}}({\mathbb C})
\oplus M_{q_n}({\mathbb C}) \hookrightarrow {\mathfrak A}_{\theta}$
with the property that $$\| \psi_n (u_{\theta_{n-1}} \oplus
u_{\theta_n}) - \psi_{n+1}(u_{\theta_n} \oplus u_{\theta_{n+1}}) \| <
\frac{2\pi}{q_{n-1}}$$ and $$\| \psi_n (v_{\theta_{n-1}} \oplus
v_{\theta_n}) - \psi_{n+1}(v_{\theta_n} \oplus v_{\theta_{n+1}}) \| <
\frac{\pi}{q_{n-1}} + \frac{4\pi}{q_n}.$$ Due to the exponential growth
of the $q_n$'s we have that both $\{ \psi_n (u_{\theta_{n-1}} \oplus
u_{\theta_n}) \}$ and $\{ \psi_n (v_{\theta_{n-1}} \oplus
v_{\theta_n})\}$ are Cauchy sequences of unitaries and hence converge
to some unitaries in ${\mathfrak A}_{\theta}$ which satisfy the same
commutation relation as $U_{\theta}$ and $V_{\theta}$.  By
universality there exists a (necessarily faithful -- by simplicity) 
$*$-homomorphism $\Psi : A_{\theta} \hookrightarrow {\mathfrak
A}_{\theta}$ with the property that $\psi_n (u_{\theta_{n-1}} \oplus
u_{\theta_n}) \to \Psi(U_{\theta})$ and $\psi_n (v_{\theta_{n-1}}
\oplus v_{\theta_n}) \to \Psi(V_{\theta})$.  More importantly, we have
the following estimates: 
\begin{eqnarray*}
\| \Psi(U_{\theta}) - \psi_n (u_{\theta_{n-1}} \oplus u_{\theta_n})\| 
& \leq &
\sum_{k = 0}^{\infty} \| \psi_{n+k} (u_{\theta_{n+k-1}} \oplus
u_{\theta_{n+k}}) - \psi_{n+k+1}(u_{\theta_{n+k}} \oplus
u_{\theta_{n+k+1}}) \|\\ 
& \leq & \sum_{k = 0}^{\infty} \frac{2\pi}{q_{n+k-1}}\\ 
& = & 2\pi( \frac{1}{q_{n-1}} + \frac{1}{q_{n}} + \sum_{k = 0}^{\infty} 
\frac{1}{q_{n+k+1}})\\ 
& \leq & 2\pi( \frac{1}{q_{n-1}} + \frac{1}{q_{n}}) + 
\frac{4\pi\sqrt{5}}{\sqrt{5} - 1}\frac{1}{q_{n+1}}. 
\end{eqnarray*}
Note that we applied Lemma \ref{thm:lemma1} to get the last
inequality.  A similar argument shows that $$\|\Psi(V_{\theta}) -
\psi_n (v_{\theta_{n-1}} \oplus v_{\theta_n})\| \leq
\frac{\pi}{q_{n-1}} + \frac{5\pi}{q_n} + \frac{10\pi\sqrt{5}}{\sqrt{5}
- 1}\frac{1}{q_{n+1}}.$$

With these estimates in hand it is immediate that 
\begin{eqnarray*}
\|\Psi(H_{\theta}) - \psi_n(h_{\theta_{n-1}} \oplus h_{\theta_n}) \|
& \leq & (|\alpha_1| + |\alpha_2|)\big(2\pi( \frac{1}{q_{n-1}}
+ \frac{1}{q_{n}}) + \frac{4\pi\sqrt{5}}{\sqrt{5} -
1}\frac{1}{q_{n+1}}\big)\\
& + & (|\beta_1| +
|\beta_2|)\big(\frac{\pi}{q_{n-1}} + \frac{5\pi}{q_n} +
\frac{10\pi\sqrt{5}}{\sqrt{5} - 1}\frac{1}{q_{n+1}}\big).
\end{eqnarray*}

The estimate above, though a bit ugly, is the one which should be used
in practice as the constants are  smaller than those claimed in
the theorem.  To get the cleaner statement of the theorem we first
note that the recursion formula $q_{n+1} = a_{n+1}q_n + q_{n-1}$
implies $$\frac{1}{q_{n+1}} \leq \frac{1}{q_{n-1}} +
\frac{1}{q_{n}},$$ and thus 
\begin{eqnarray*}
\|\Psi(H_{\theta}) - \psi_n(h_{\theta_{n-1}} \oplus h_{\theta_n}) \| &
\leq & 2M\bigg( 2\pi + \frac{4\pi\sqrt{5}}{\sqrt{5} - 1}\bigg)(
\frac{1}{q_{n-1}} + \frac{1}{q_{n}})\\ 
& + & 2M\bigg( 5\pi + \frac{10\pi\sqrt{5}}{\sqrt{5} - 1}\bigg)(
\frac{1}{q_{n-1}} + \frac{1}{q_{n}})\\
& = & 2M\bigg(7\pi + \frac{14\pi\sqrt{5}}{\sqrt{5} - 1}\bigg)(
\frac{1}{q_{n-1}} + \frac{1}{q_{n}})\\
& = & 14\pi M\frac{3\sqrt{5} - 1}{\sqrt{5} - 1}(
\frac{1}{q_{n-1}} + \frac{1}{q_{n}})\\
& \leq & 204M(\frac{1}{q_{n-1}} + \frac{1}{q_{n}}).
\end{eqnarray*}

Now we quote some general results concerning approximation
of spectral quantities.  Given two operators $S, T \in B(H)$ one has: 
\begin{enumerate}
\item (cf.\ \cite[Theorem 3.27]{HRS}) For every $\varepsilon > 0$,
$\sigma^{(\varepsilon)}(S) \subset \sigma^{(\varepsilon+\|S - T\|)}(T)
\subset \sigma^{(\varepsilon+2\|S - T\|)}(S)$.


\item If both $S$ and $T$ happen to be normal then $d_H (\sigma (S),
\sigma (T)) \leq \|S - T\|$.
\end{enumerate}
With these general facts and the norm estimates above one easily
completes the proof by observing that the spectral quantities of $S
\oplus T$ (i.e.\ pseudo or actual spectrum)
are just the union of the spectral quantities for $S$ and $T$.
\end{proof}

\begin{rem} {\em Almost tri-diagonality.} Note that $h_{\theta_n}$ is
``almost'' a tri-diagonal matrix.  It is non-zero in the upper-right
and lower-left entries but otherwise is tri-diagonal.  In particular,
these matrix models are sparse.
\end{rem}

\begin{rem} {\em Discretized Schr\"{o}dinger operators with polynomial
potential. }

The theorem above is stated for a class of operators that includes all
of the examples discussed in \cite{boca:book} (i.e.\ Almost Mathieu
operators and a few non-self-adjoint examples).  However, similar
results hold for a broader class of examples that includes all the
one-dimensional discretized Schr\"{o}dinger operators with polynomial
potential.  In fact, if $R(X, Y)$ is any polynomial in non-commuting
variables $X$ and $Y$ then $$\| \Psi( R(U_{\theta}, V_{\theta})) -
\psi_n( R(u_{\theta_{n-1}}, v_{\theta_{n-1}}) \oplus R(u_{\theta_{n}},
v_{\theta_{n}}) )\| = O(\frac{1}{q_{n-1}} + \frac{1}{q_{n}}).$$
Unfortunately, the constants involved get quite complicated quite
quickly but they can be explicitly worked out for a given polynomial
should one so desire.  Having this norm estimate one gets similar
convergence results for spectral quantities in this setting.  Finally
we should point out that the spectral quantities of {\em any} operator
in an irrational rotation algebra (e.g.\ Schr\"{o}dinger operators
with arbitrary quasiperiodic potential) can, in principle, be
approximated with this method -- first approximate by a polynomial and
then proceed -- however, it seems impossible to control rates of
convergence in this generality.
\end{rem}

\begin{rem} {\em These rates of convergence are sharp.}

To see this we consider the case that $H_{\theta} = U_{\theta}$ and
$\theta$ is an irrational number for which the $a_n$'s are uniformly
bounded.  Then $d_H (\sigma(H_{\theta}), \sigma(h_{\theta_{n-1}}) \cup
\sigma(h_{\theta_n})) \leq 292(\frac{1}{q_{n-1}} + \frac{1}{q_{n}})$.
It is easy to check that $\frac{1}{q_{n-1}} + \frac{1}{q_{n}} \leq
\frac{4 + a_n}{q_{n-1} + q_n}$ and hence if one could prove that $d_H
(\sigma(H_{\theta}), \sigma(h_{\theta_{n-1}}) \cup
\sigma(h_{\theta_n})) = o(\frac{1}{q_{n-1}} + \frac{1}{q_{n}})$ then
it would follow that $d_H (\sigma(H_{\theta}),
\sigma(h_{\theta_{n-1}}) \cup \sigma(h_{\theta_n})) =
o(\frac{1}{q_{n-1} + q_{n}})$.  However this is impossible as it would
imply that the Lebesgue measure of $\sigma(H_{\theta})$ is zero (which
it isn't: $\sigma(U_{\theta}) = \{ \lambda \in {\mathbb C}: |\lambda|
= 1\}$) since there are at most $q_{n-1} + q_n$ points in
$\sigma(H_{\theta_{n-1}}) \cup \sigma(H_{\theta_n})$. 
\end{rem}

\begin{rem} {\em Implementability.}
It is, of course, true that for some irrational numbers
the strategy proposed in this paper for approximating spectra can't
reasonably be implemented.  For example, it can happen that $q_1$ is very
small while $q_2$ is far larger than any computer can handle (just
take $a_2$ to be huge).  On the other hand, the set of irrationals for
which the associated sequences $\{a_n\}$ stay bounded is a dense set
in $(0,1)$.  For example, all irrationals which are the roots of
quadratic equations have this property; hence $\sqrt{r}$ for any
rational number and these are already dense in $(0,1)$.  In other
words, for a dense set of irrationals this program is reasonable and
hence the question becomes how spectra behave as the parameter
$\theta$ changes.  But it is well known that spectra of operators in
irrational rotation algebras vary continuously in $\theta$ (see the
next section) and hence the program suggested in this paper is about
as practical as could be hoped for.
\end{rem}

\begin{rem} {\em Uniqueness of convergents.}  Another basic, but
remarkable, property of continued fractions is that if $\theta$ is
irrational and $| \theta - \frac{p}{q}| < \frac{1}{2q^2}$ then
$\frac{p}{q}$ is necessarily one of the convergents in the continued
fraction decomposition of $\theta$ (i.e.\ if $\theta = [a_1, a_2,
\ldots]$ then there is an $n$ such that $\frac{p}{q} = [a_1, \ldots,
a_n] = \frac{p_n}{q_n}$ -- see \cite[Theorem 184]{hardy-wright}).  The
point is that if one {\em starts} with a rational $\frac{p}{q}$ and
only looks at the irrationals $\theta$ which are a distance at most
$\frac{1}{2q^2}$ away from $\frac{p}{q}$ then the spectral quantities
of $h_{\frac{p}{q}}$ will provide good approximations to those of
$H_{\theta}$ (since $\frac{p}{q}$ necessarily arises in the continued
fraction decomposition of all such $\theta$ and hence the theorem
above applies).  Of course, rates of convergence become trickier in
this setting but our only point is that one can simultaneously
approximate a small interval of operators with a single matrix model.
\end{rem}

\section{General One-sided Rates of Convergence}

In our final section we will record a consequence of the following
theorem of Haagerup and R{\o}rdam (cf.\ \cite[Corollary
3.5]{boca:book}).  Among other things this shows that if one fixes a
polynomial in the canonical generators then the spectral quantities
vary continuously in $\theta$.    

\begin{thm} For every $\theta, \theta^{\prime} \in [0,1)$ there exists
a Hilbert space $H$ and injective $*$-homomorphisms $\pi_{\theta} :
A_{\theta} \to B(H)$, $\pi_{\theta^{\prime}} : A_{\theta^{\prime}} \to
B(H)$ such that $$\| \pi_{\theta^{\prime}}(U_{\theta^{\prime}}) -
\pi_{\theta}(U_{\theta})\|, \|
\pi_{\theta^{\prime}}(V_{\theta^{\prime}}) -
\pi_{\theta}(V_{\theta})\| \leq 9\sqrt{6\pi|\theta-
\theta^{\prime}|}.$$
\end{thm}

As in the last section we will restrict attention to operators of the
form $H_{\theta} = \alpha_1 U_{\theta} + \alpha_{-1} U_{\theta}^* +
\beta_1 V_{\theta} + \beta_{-1} V_{\theta}^*$ as the estimates are
cleaner in this setting.  However, the same techniques handle more
general polynomials as well.  In the theorem below we only get
``one-sided'' rates of convergence (as opposed to estimates on the
Hausdorff distance).  The modest advantage, however, is that this
strategy is implementable for any irrational number.

We find the following notation convenient: For compacts $\Lambda,
\Sigma \subset {\mathbb C}$ and $\delta > 0$ we write $\Lambda
\subset^{\delta} \Sigma$ if for each $\lambda \in \Lambda$ there
exists $\sigma \in \Sigma$ such that $| \lambda - \sigma | < \delta$.
Note that $d_{H}(\Lambda, \Sigma) < \delta$ if and only if $\Lambda
\subset^{\delta} \Sigma$ and $\Sigma \subset^{\delta} \Lambda$.

\begin{thm} Let $\theta \in (0,1)$ be irrational and for each $n \in
{\mathbb N}$ choose an integer $p_n \in \{ 0,1,\ldots,n-1\}$ such that
$|\theta - \frac{p_n}{n}| \leq \frac{1}{2n}$.  If $H_{\theta}$ is an
operator polynomial as above, $h_n = \alpha_1 u_{\frac{p_n}{n}} +
\alpha_{-1} u_{\frac{p_n}{n}}^* + \beta_1 v_{\frac{p_n}{n}} +
\beta_{-1} v_{\frac{p_n}{n}}^*$ is the associated matrix model and $M
= \max\{|\alpha_{\pm 1}|, |\beta_{\pm 1}|\}$ then the following
statements hold:
\begin{enumerate}
\item For every $\varepsilon > 0$,
$$\sigma^{(\varepsilon)}(h_n) \subset \sigma^{(\varepsilon +
\frac{C_1}{\sqrt{n}})}(H_{\theta}),$$ where $C_1 \leq 36M\sqrt{3\pi}.$


\item If $H_{\theta}$ happens to be normal then
$$\sigma(h_n) \subset^{\frac{C_1}{\sqrt{n}}} \sigma(H_{\theta}),$$
where $C_1 \leq 36M\sqrt{3\pi}.$
\end{enumerate}
\end{thm}

\begin{proof} To prove this result we first observe that all of the spectral
quantities associated to $h_n$ are contained in the corresponding
quantities for $H_{\frac{p_n}{n}} = \alpha_1 U_{\frac{p_n}{n}} +
\alpha_{-1} U_{\frac{p_n}{n}}^* + \beta_1 V_{\frac{p_n}{n}} +
\beta_{-1} V_{\frac{p_n}{n}}^*$ (since universality gives us a
$*$-homomorphism taking $H_{\frac{p_n}{n}} \mapsto h_n$). The
Lip$^{1/2}$-continuity theorem of Haagerup-R{\o}rdam then gives upper
bounds on the distance between the spectral quantities of $H_{\theta}$
and $H_{\frac{p_n}{n}}$ which completes the proof.
\end{proof}

\begin{rem} {\em Hausdorff Convergence.}  
It can be shown that $d_H(\sigma^{(\varepsilon)}(h_n),
\sigma^{(\varepsilon)}(H_{\theta})) \to 0$ (for every $\varepsilon >
0$).  The proof is similar to the proof of \cite[Theorem
3.5]{brown:finitesection} and hence will be omitted.  However, it does
not seem possible to control the rate of convergence in the Hausdorff
metric.  For example, it is certainly not the case that rates of
convergence are controlled by $|\theta - \frac{p_n}{n}|^{\sigma}$ for
any $\sigma > 0$ because many irrational numbers (e.g.\ transcendental
numbers) have extremely good rational approximations.  To be more
precise, it can happen that a particular $\theta$ has the property
that there are infinitely many co-prime solutions to the inequality
$|\theta - \frac{p}{q} | < \frac{1}{q^k}$ where $k$ is some fixed
natural number.  If the Hausdorff distance between spectra could be
bounded above by $|\theta - \frac{p}{q} |^{\sigma}$ then we would also
get an upper bound of $(\frac{1}{q})^{k\sigma}$.  However, when
$k\sigma > 1$ it would then follow that the spectrum of $H_{\theta}$
has Lebesgue measure zero since there are at most $q$ eigenvalues for
the matrix $h_{\frac{p}{q}}$.  Since these operators need not have spectra of
measure zero it follows that an upper bound of the form $|\theta -
\frac{p_n}{n}|^{\sigma}$ is impossible.
\end{rem}

\bibliographystyle{amsplain}

\providecommand{\bysame}{\leavevmode\hbox to3em{\hrulefill}\thinspace}

\end{document}